\documentclass [a4paper,14pt ]{extarticle}
\usepackage[dvips,pdftex]{graphicx,graphics}
\usepackage{amsfonts}
\usepackage[utf8]{inputenc}
\begin{document}
 \baselineskip=16pt
\clearpage

\newcommand{\N}{\mathbb{N}}
\newcommand{\qq}{\qquad}
\newcommand{\q}{\quad}

\newcommand{\si}{\sigma}
\newcommand{\la}{\lambda}
\newcommand{\al}{\alpha}
\newcommand{\be}{{\beta}}
\newcommand{\ve}{{\varepsilon}}

\baselineskip=16pt

 \setcounter{page}{1}

{\ }

\vspace{3ex} 

 MSC-2010-class: 11A25, 11N56 (Primary); 11M26 (Secondary)

\begin{center}

\vspace{3ex}

{\ } {\bf  \Large Remainder in the Modified Mertens Formula

\vspace{1ex}
 and Ramanujan Inequality\footnote{\ \ \ This work
 was supported
by the grant of Russian Foundation of Fundamental Research\\
$\qq {\ } {\qq } {\ } $ (project \# $14-01-00684.$)}
}

\vspace{2ex}

{\bf \large by Gennadiy Kalyabin\footnote{\q Samara, 
 Russia; \ gennadiy.kalyabin@gmail.com}}

\end{center}
\vspace{1ex} 

{\it Abstract:} 
A highly strong upper estimate in the modified asymptotic 

formula for sums of the primes' reciprocals is proved   to be 
necessary

 (as well as  sufficient) in order the Ramanujan inequality holds true.  
 
 Some other criteria in similar terms are also obtained.

\vspace{1ex}


 inequality, Riemann Hypothesis

\vspace{1ex}

{\it Bibliography}: 9 items

\vspace{3ex}

{\bf 1. \ Notations, brief history  and main results}

\vspace{2ex}

 As usually, let $\N$ be a set  of all positive integers, 
 $\N_0:=\N \cup\{0 \}$, $p$ run the set 
  $\mathbb{P}:=\{p_1, p_2, \dots \},\ p_j<p_{j+1},$ of all primes, 
$\ve$ is an arbitrary positive number, $C_y$ stand for positive constants 
which may depend only on a parameter $y$; \ symbols $\triangleright$ and $\Box$
denote the proof's beginning and end; \
 $\log x$ and $\gamma$ stand (resp.) for the natural logarithm of a positive $x$ and the Euler-Masceroni constant:
\vspace{-2ex}
  $$ \gamma:= \lim_{n\to\infty} \left( \sum_{k=1}^n \frac{1}{k}
   - \log n\right)=0.577\ 215\ 664\ \dots
  \eqno{(1.1)}
  $$
\vspace{-2ex}  

  In 1874 F. Mertens [1] proved his famous asymptotic formula
  \vspace{-2ex}  
  $$  S(x) :=\sum_{p\le x} \log \frac{p}{p-1}
  = \log \log x + \gamma+ R(x) \hbox{ with  }
   R(x)= O\left(\frac{1}{\log x}\right). 
 \eqno{(1.2)}
  $$
 
\vspace{-1ex}
 
 The best known  {\it unconditional}, (i. e. without assumption of 
the Riemann Hypothesis {\bf (RH)} ), estimate for this remainder 
at the moment (2021) seems to be   
 $R(x)=O(\exp(- c(\log x)^{3/5}(\log \log x)^{-1/5})$.
  
\vspace{1ex}

  In 1984 {\it assuming } {\bf RH}
   G.~Robin [2, Th. 3] has come to the fundamentally stronger
    estimate: $|R(x)|<\log x/(8\pi\sqrt{x}), \ x>X_0$.

  \vspace*{\fill}
  \clearpage
   
 \vspace{1ex} 
  
 We will present an integer $N>1$  as its canonical 
  factorization in primes
  $$ N:=p_1^{\al_1} p_2^{\al_2} \dots p_k^{\al_k}; 
  \qq \al_j\in\N_0,\ \al_k>0,
  \eqno{(1.3)}
  $$
  where the number $k:=k(N)$ and exponents 
  $\al_j:=\al_j(N),\ 1\le j \le k,$ are uniquely defined by $N$. 
 
 \vspace{1ex} 
 
  The greatest prime factor $p_k:=p_{k(N)}$ of $N$ 
  will be denoted by $\hbox{\bf gpf}(N)$.  
  
 \vspace{1ex}

  Let $\sigma(N)$ stand for the arithmetic multiplicative function sum of all divisors of $N\in \mathbb{N}$. The  properties of this function are well described in a very informative paper [3], which contains a lot 
  of valuable historical remarks, as well as many definions, notations 
  and facts widely used in this paper. In particular, in Sect. 5   the classical formula  for $\si(N)$  is adduced, $N$ being defined in (1.3), 
  namely:
 \vspace{-1ex}
$$ \si(N):= \prod_{j=1}^{k} (1+p_j +\dots +p_j^{\al_j}) 
= \prod_{j=1}^{k}  \frac{p_j^{\al_j+1}-1}{p_j -1}
  \eqno{(1.4)}
  $$

      T. Gronwall in 1913, basing on (1.2) established the  sharp  upper 
      order of $\si(n)$, namely he proved  [4] that: 
      \vspace{-2ex} 
     $$ \limsup_{N\to\infty} G(N)  = e^{\gamma} = 1.781\ 072\dots; \hbox{   where   }    G(N):= \frac{\si(N)}{N \log \log N},
  \eqno{(1.5)}
  $$
 which we will  call  Gronwall numbers.
\vspace{1ex}

S. Ramanujan has noticed (in 1915, 
the first publication  in 1997 [5]) that:
\vspace{1ex}

{\it if {\bf RH}  holds true, then in addition to (1.6) for all $N$ sufficiently large the following {\it strict} (Ramanujan) inequality 
{\bf (RI)} takes place:}  
$$ G(N)<e^{\gamma}, \qq \forall N>n_0.
\eqno{(1.6)}
$$

 \vspace{-1ex}
 Almost 70 years later G. Robin [2, Th 1] proved a  paramount  assertion, which in a sense complements the Ramanujan's result, namely:\  

\vspace{1ex}
{\it if  (1.6) holds true for all integers  $N>5040$, then {\bf RH} is valid. } 

\vspace{1ex}

 We will call (1.6) with $n_0=5040$ the {\it Ramanujan-Robin} inequality  {\bf (RRI)}, in which the statement  of {\bf RH} is exhaustively  encoded in terms of $\si(N)$. 
 
Robin has also shown that {\it in fact {\bf RI} (with uncertain $n_0$)  and  {\bf RRI}   are equivalent}, because  if {\bf RH} holds false, then there are infinitely many $N$'s such that $G(N)> \gamma$.

 \vspace{1ex} 
 
In this paper we {\it do not strive} to prove any of these conjectures but
rather to reveal the direct interrelation between them and the remainder in the {\it modified} Mertens formula:
$$  S(x)  = \log \log \theta(x) + \gamma+Q(x),
 \eqno{(1.7)}
  $$
   which differs from (1.2) by replacing $x$ in $\log\log$ by    the first Chebyshev function $\theta(x):=\sum \{\log p: p\le x\} $ (cf [6, 3.1]).

\vspace{3ex}
  
\vspace*{\fill}
\clearpage

{\bf Theorem.}\  {\it {\bf RRI} 
 is equivalent to each of the following three  conditions:}
\vspace{-1ex}
  $$   \forall \varepsilon>0\ \forall x>1:  \
 Q(x)  < C_{\varepsilon}x^{-0.5 + \varepsilon} \
 \eqno{(1.8)}
 $$
 $$
\  \forall \varepsilon>0\ \forall x>1:  
  \ Q(x) > -  C_{\varepsilon}x^{-0.5 + \varepsilon},
 \eqno{(1.9)}
  $$ 
\vspace{-1ex} 
\vspace{-1ex}
$$ \q
 A_0 := \limsup_{ x\to +\infty}\ Q(x) \sqrt{x} \log x  < +\infty.  
  \eqno{(1.10)}
  $$ 
  
 { \it In addition,  (1.10) necessarily implies that 
 $A_0\le 2\sqrt{2}$; for this reason the situation
  $2\sqrt{2}<A_0<+\infty$ is {\it logically} impossible. }

\vspace{1ex}
  
The proof of the Theorem is set forth  in Sect. 3; in Section 2 all needed auxillaries and the main  Lemma are adduced; in Sect. 4 
some corollaries and directions of further research are given.

\vspace{3ex}

{\bf 2. \q \ Some known facts and main lemma}

\vspace{2ex}

Further the well-known  assertions are brought together concerning the asymptotic 
behavior of primes  [6, Ch. 5] in their weak form  sufficient for our purposes:

\vspace{1ex}

{\bf Proposition 1. (i)} {\it For  all  $x>1$ one has}
  $ |\theta (x) -x| <C_0x/\log x;$ 

\vspace{1ex}
 
{\bf (ii)} $  \ p_{k+1} - p_k < C_1p_k/ \log p_k$; 

\vspace{1ex}

{\bf (iii)} {\bf RH} {\it is equivalent to each of the two relationships}:
\vspace{-2ex} 
$$ \forall \ve > 0 \exists C_{\ve} \forall x>1:
 |\theta(x) -x| < C_{\ve}x^{0.5+\ve}; \
\ |\theta(x) - x| <\frac{ \sqrt{x} \log^2x}{8\pi},\ x>X_0.
\eqno{(2.1)} 
 $$

For the sequel we will need the  (perhaps also well-known) ascertion, 
which follows from Proposition 1{\bf(i)}: 

\vspace{1ex}

{\bf Proposition 2. } {\it Let $\la>1$; then for all 
$x> X_{\la}:=\exp(\max(1,2/(\la-1)))$ one has:}
$$ Y=Y(x,\la):= \sum_{p>x} \frac{1}{p^{\,\la}}
= \frac{1+\delta(x, \la)}{(\la-1) x^{\la-1}\log x}; 
\q |\delta(x, \la)|<\frac{C_1}{(\la-1)\log x}.
\eqno{(2.2)} 
 $$

\vspace{1ex}
$\triangleright$ In fact,  using the integration by parts one obtains
 $$ Y=\int_{x^+}^{\infty} \frac{d\theta (t)}{t^{\la}\log t}
= \frac{\theta(t)}{t^{\la}\log t}\ \Bigg|_{x^+}^{\infty} 
- \int_x^{\infty} \theta(t) \left(\frac{1}{t^{\la}\log t}\right)^{\prime} dt                                                                                                                                                                                                                                                                                                                          
$$
$$ =-\frac{\theta(x^+)}{x^{\la}\log x}
 +\int_x^{\infty}\frac{\theta(t) (\la\log t +1)}{t^{\la+1} \log^2t}\, dt 
\eqno{(2.3)} 
 $$ 

\vspace*{\fill}
\clearpage

Analogously, replacing here $\theta(t)$ by $t$ one obtains the identity:
 $$ J=J(x,\la):=\int_{x}^{\infty} \frac{dt}{t^{\la}\log t}
 =-\frac{x}{x^{\la}\log x}
 +\int_x^{\infty}\frac{t\ (\la\log t +1)}{t^{\la+1} \log^2t}\, dt 
$$
$$ = -\frac{1}{x^{\la-1}\log x} +\la J +\frac{\be J}{\log x};
 \hbox{  where  } 0<\be=\be(x, \la)<1\hbox{  for  }  x>X_{\la},
\eqno{(2.4)} 
 $$ 
whence it follows that for all $x>X_{\la}$ (explanations below):
 $$ J(x, \la) =  \frac{1}{ (\la -1) x^{\la-1} \log x
 \left(1-\frac{\be}{(\la-1) \log x } \right)}
$$
$$ \Rightarrow\q 0<J(x, \la) -  \frac{1}{ (\la -1) x^{\la-1} \log x} < 
\frac{2}{ (\la -1)^2 x^{\la-1} \log^2 x}.
\eqno{(2.5)} 
 $$ 

Here we have taken into accout that since $x>X_{\la}$ then by  virtue of 
the $X_{\la}$-definition the number $t:=1/(\la-1)\log x $ belongs to the interval  $(0, 1/2)$ and hence the inequality $1/(1-t)<1+2t$ holds.

\vspace{1ex}

On the other hand, substracting (2.4) from (2.3) and using 
Proposition 1{\bf(i)} leads to:
 $$ |Y-J| \le \frac{|\theta(x^+) - x|}{x^{\la}\log x} 
+\int_x^{\infty}\frac{|\theta(t) -t| (\la\log t +1)}{t^{\la+1} \log^2t}\, dt
$$
$$ \le \frac{C_0}{\log x}\left( \frac{2}{x^{\la-1}\log x}+  J \right)
< \frac{4C_0}{(\la-1) x^{\la-1}\log^2 x}, \q \forall x>X_{\la}.
\eqno{(2.6)} 
 $$

\vspace{1ex}

Joining (2.6) with (2.5) one comes to (2.2) $\Box$.

\vspace{2ex}

The main role in the proof of the Theorem plays the following 
{ unconditional} assertion, binding Mertens function $S(x)$ 
and Gronwall numbers $G(N)$, which is the most important 
and complicated part of the  paper.

\vspace{1ex}

{\bf Lemma.} {\it For any  $k\in \N$ there are a real number 
$\delta_k, \ (\{\delta_k\} \to 0$ as $ k\to\infty),$ and an integer 
$N_k^*$ such  that} $ \hbox{\bf gpf}(N_k^*)=p_k $ {\it and}
\vspace{-1ex}
  $$ \log G(N_k^*) > S(p_k) - \log\log \theta(p_k) 
  - \frac{2\sqrt{2}+\delta_k }{\sqrt{p_k}\log p_k}. 
\eqno{(2.7)} 
 $$

$\triangleright$ 1) First we { describe} the  special construction of
$\{N_k^*\}_{k=1}^{\infty}$ providing (2.7).
 
\vspace{1ex}

We'll suppose that $k$ is large enough; put $ r=r_k:=[\sqrt{ \log 2 p_k }],  $ and define:
\vspace{-1ex}
$$   \ \ q_1:=p_k, \q q_m:= \max\{p_j: p_j^m \le {2p_k}\}, \ 2\le m\le r. 
\eqno{(2.8)} 
 $$ 

\vspace{3ex}
 
\vspace*{\fill}
\clearpage

\vspace{-1ex}

 In other words, $q_m =q_{m,k}$ is the greatest prime 
  $\le(2p_k)^{1/m}$; hence   $q_{m-1}< q_{m}$ for all $m, 1< m\le r.$ 
From Proposition 3 one may easily deduce  that the quantity 
$q_{m,k}=(2p_k)^{1/m}(1-\delta_{k,m});
\ 0\le \delta_k :=\max_{1<m\le r} \delta_{k,m} \to 0,\ k\to\infty$.

Let $\nu=\nu_r:=\max\{j: p_j\le q_r\}, \ H:=q_r^{r+1}$; define the exponents $\{\al_j\}_{j=1}^k$
$$ \al_j := \left[\frac{\log H}{\log p_j}\right] -1,\ \ \hbox{if  } j<\nu;
\q \al_j:= \max \{m\le r: q_m\ge p_j\}, \ \hbox{if  } j\ge\nu; 
\eqno{(2.9)} 
 $$ 

It is clear that: 1) $\al_{\nu}=\al_{\nu+1}=r,\ p_{\nu}=q_r,$
\ 2) $\al_j\ge \al_{j+1},\ 1\le j<k$,  \\
3) the equality $\al_j =m<r$ is equivalent to   $q_{m+1} < p_j \le q_m.$

\vspace{1ex}

Let $ T(x):=\exp(\theta(x))$ stand for a product of all primes 
$p\le x.$
 \vspace{1ex} 

Now we are able to determine the numbers $N_k^*,$ for which  
the relationship (2.7) is guaranteed:
 \vspace{-2ex} 
  $$ N_k^*: =\prod_{m=1}^{r} T(q_m) \cdot
 \prod_{j=1}^{\nu-1} p_j^{\al_j-r}
  = \prod_{j=1}^k p_j^{\al_j}. 
\eqno{(2.10)} 
 $$ 

{2) } Let's study  the quantity 
 $\eta=\eta_k:=\log N_k^*=E_k+F_k; E_k:= \sum_{m=1}^r \theta(q_m),
\\  F_k:=\sum_{j=1}^{\nu-1} (\al_j-r) \log p_j.$ 
Having taken into account the definition (2.8) of $q_m$ and 
 the relationships:\ 
  $ \max\{ |1-\theta(q_{m,k})(2p_k)^{-1/m}|: 1<m\le r\} \to0,\ k\to\infty,$
$\nu<p_{\nu}=q_r<(2p_k)^{1/r}, \log H=(r+1)\log q_r,$  
 one has:
\vspace{-1ex} 
 $$  \ E_k:=  \theta(p_k)+C_k\sqrt{p_k} 
+ O(p_k^{1/3} \sqrt{\log 2p_k}),\ C_k\to\sqrt{2}, \ k\to\infty;
$$
$$\qq 0<F_k < \nu\log H < q_r (r+1)  \log q_r  =O(p_k^{\ve})
\ \Rightarrow \eta_k -\theta(p_k) \approx \sqrt{2p_k}. 
\eqno{(2.11)}
  $$

 3)  From (1.4) and (1.2) it follows  that
  \vspace{-1ex} 
    $$ \log G(N_k) =  \log \frac{\si(N_k)}{N_k} -\log\log\log N_k= \sum_{j=1}^k \log\ \frac{p_j^{\al_j+1}-1}{ p_j^{\al_j}(p_j-1)}  
  - \log\log \eta_k
  $$
  \vspace{-1ex}
$$= \sum_{j=1}^k \log\  \frac{p_j}{p_j-1}  
  - \sum_{j=1}^k \log \frac{p_j^{\al_j+1}}{p_j^{\al_j+1}-1}   
   - \log\log \eta_k = S_k - U_k - V_k.
   \eqno{(2.12)} 
    $$

 Now  with certain  $t_k$ in between of $\theta(p_k)$ and $\eta_k$,
 one has: 
\vspace{-1ex}
$$V_k - \log\log \theta(p_k)
=\frac{\eta_k-\theta(p_k)}{t_k\log t_k} 
\approx \frac{\sqrt{2}}{ \sqrt{p_k}\log p_k},\q k\to\infty.
\eqno{(2.13)} 
 $$ 

\vspace{1ex}

4) To make sure that  the quantity $U_k$ is also 
$\approx\sqrt{2}/\sqrt{p_k}\log p_k$ as $k\to\infty$,
 we   present it as a sum: 

\vspace*{\fill}
\clearpage

\vspace{-2ex}
 $$ U_k =
\sum_{m=1}^r U_{k, m};
 \qq  U_{k, m}:=\sum_{q_{m+1}<p_j\le q_{m}} \log\frac{p_j^{m+1}}{p_j^{m+1}-1}, \ m<r; 
$$
$$\q U_{k, r}:=\sum_{j=1}^{\nu-1} \log\frac{p_j^{\al_j+1}}{p_j^{\al_j+1}-1}.
\eqno{(2.14)} 
$$

\vspace{1ex}

All summands $U_{k, m}$ here are positive. Using the elementary inequality:  $-t^2<\log (1-t)+t<0, \ 0< t < 1/4$,   easily deduced from the Taylor formula, one may assert   that for $m<r,\ k>k_0$  and some 
$\delta_{j,m}\in (0, 1)$:
$$ U_{k, m}= 
\sum_{q_{m+1}<p_j \le q_{m}}  -\log \left(1-\frac{1}{p_j^{m+1}}\right)
 = 
\sum_{q_{m+1}<p_j\le q_m}
\left(\frac{1}{p_j^{m+1}} + \frac{\delta_{j,m}}{p_j^{2m+2}}\right),
\eqno{(2.15)}
$$

\vspace{1ex}

5)  Recollecting now  the definition (2.2)  of the quantity $Y(x, \la)$ 
 in Proposition 2, we may rewrite the latter equality as follows:
$$ U_{k,m}=Y(q_{m+1}, m+1) - Y(q_{m}, m+1)+W_{k, m}; 
$$
$$
\q 0<W_{k, m}<Y(q_{m+1}, 2m+2), \ 1\le m<r.
\eqno{(2.16)}
$$

Applying  the relationship (2.2) with $\la=m+1,\ 2m+2,\ x=q_{m+1}, 
q_{m}, $ and having taken into account that $q_{m+1}\approx(2p_k)^{1/(m+1)}$, by virtue of defining formula (2.8),   one comes to the estimates 
\vspace{-1ex}
$$ U_{k,m} < \frac{1}{mq^{m}_{m+1} \log q_{m+1}}
\left(1+\frac{C_1}{\log  q_{m+1}}\right)<{C_2}\ p_k^{-m /(m+1)}; 
\q m<r.
\eqno{(2.17)}
$$

 Further, for $m=r$  due to the fact that $\al_j\log p_j > (r+1)\log q_r-\log p_{\nu}$ for $j<\nu$ (cf. the left  part of definition (2.9)) and $\nu<p_{\nu}=q_r$,   one obtains:
\vspace{-1ex}
 $$  U_{k,r}<{2} \sum_{j=1}^{\nu -1}  \frac{1}{p_j^{\al_j+1}}
<\frac{2\nu p_{\nu}}{q_r^{r+1}} < 
 \frac{2}{q_r^{r-1}}<\frac{3}{(2p_k)^{1-1/r}} = O(p_k^{-1+\ve}). 
\eqno{(2.18)}
$$

From these two estimates it follows that
for $k$ large enough:
\vspace{-1ex}
$$ \sum_{m=2}^{r} U_{k,m} < C_2 p_k^{-2/3} (\log p_k)^{1/2}.
\eqno{(2.19)}
$$

\vspace{-2ex} 
6) At last, if $m=1$, then again by virtue 
of Proposition 2, one has 
\vspace{-1ex}
$$ U_{k,1}=Y(q_2, 2) +O(p_k^{-1})
=\frac{1+O(1/\log x)}{\sqrt{2p_k}\log\sqrt{2p_k}}
=\frac{\sqrt{2}+O(1/\log p_k)}{\sqrt{p_k}\log p_k}.
\eqno{(2.20)}
$$
\vspace*{\fill}
\clearpage

\noindent
whence in junction with  (2.13) and (2.19) it follows that
 $U_k\approx\sqrt{2}/\sqrt{p_k}\log p_k$, and joining this with (2.11), 
(2.12), one comes to the limit relationship:
\vspace{-1ex}
$$ \log G(N_k^*) - (S(p_k) - \log\log \theta(p_k))\approx
 \frac{2\sqrt{2}}{\sqrt{p_k}\log p_k},\q k\to \infty,
\eqno{(2.21)}
$$
 which in turn implies (2.7) $\ \Box$. 

\vspace{1ex}

Now we have got all the tools needed  to move forward.

\vspace{2ex}

{\bf 3. \q \ Proof of the Theorem}

\vspace{2ex}

{\bf Sufficiency.\ } 
$\triangleright$\ Due to  Nicolas result (cf [7], [2], Sect. 4) the
 {\it negation} of   {\bf RH} implies 
   $Q(x) = \Omega_{\pm}(x^{-b})$ for some    $b\in (0, 0.5)$, 
   i. e. according to the meaning of the symbol          
             $\Omega_{\pm}$,  for some $\delta>0$ 
  and any $X>0$ there are  $y>z>X$ such that
                     $Q(y)>\delta y^{-b}, \  Q(z)<-\delta z^{-b}$, but
each of these two inequalities contradicts (resp.) to (1.8), (1.9). 
Besides, obviously (1.10)\ $\Rightarrow$ (1.8).

 Thus it is proved that each of (1.8), (1.9) and (1.10)
 implies {\bf RH}  $\Box$.

\vspace{1ex}

{\bf Necessity.\ } $\triangleright$ Let us {\it suppose}  that (1.10) is false,
or more precisely, that $B_0 > 2\sqrt{2}$;
 then taking into account 
the relationships   (1.2), (1.7) and (1.10), one may conclude that for any
fixed   $\varepsilon_1\in (0, B_0 - 2\sqrt{2})$  the set   
 \vspace{-1ex}
$$ K_{\varepsilon_1}:=\left\{ k: Q(p_k) > 
 \frac{2\sqrt{2} + \varepsilon_1}{\sqrt{p_k}\log p_k} \right\} 
 \q \hbox{is {\it infinite.} } 
 \eqno{(3.1)} 
 $$

But then  by virtue of Lemma and the equality
 $Q(x)=S(x) -\log\log \theta(x) - \gamma,$  (cf (1.7), (2.1)) 
one obtains for all sufficiently large  $k\in {K}_{\varepsilon_1}$ 
 $$ \log G(N_k^*) > \gamma + Q(p_k) - \frac{2\sqrt{2} +\delta_k }{\sqrt{p_k}\log p_k} 
 >\gamma+  \frac{\varepsilon_1 - \delta_k}{\sqrt{p_k}\log p_k} 
 >\gamma,
\eqno{(3.2)} 
 $$ 
 because $\delta_k\to0,\ k\to\infty,$ whereas $\varepsilon_1>0,$ and consequently {\bf RH} holds false. 

\vspace{1ex}

Further, assuming {\ RH} one deduces from (1.2), (1.7)
 and Proposition 1{\bf(iii)}, that there is $t, (t-x)(t-\theta(x))<0$ for which  
$$ |R(x) - Q(x)| = \frac{|\theta(x) -x|}{t\log t} <
 \frac{1+o(1)}{8\pi \sqrt{x} \log  x},
\eqno{(3.3)} 
 $$ 
and joining this with Robin's estimate 
$|R(x)|< \log x /(8\pi\sqrt{x})$, mentioned in Sect. 1, 
one obtains $|Q(x)|< (1+o(1)) \log x /(8\pi\sqrt{x})$, 
which in turn implies (1.9)  \ $\Box$.

\vspace{1ex}

This completes the Theorem's proof.

\vspace*{\fill}
\clearpage

{\bf 4. \ Conclusive Remarks.}

\vspace{1ex}

1) The assertions  in Theorem may be presented in a {\it discrete}
 form, when $x$ in (1.2) runs only  the sequence of primes 
$\{p_j\}_{j=1}^{\infty}$.

\vspace{1ex}

2) One may also  replace  in  (1.2) $\log(p/p-1)$ by primes reciprocals $1/p$.

\vspace{1ex}

{\bf Corollary.}
 {\it  The Ramanujan  inequality (1.6) (and thus {\bf RH})  
  is equavalent to each of the following three unilateral  estimates for all  suficiently large $k$: }

\vspace{-2ex}

 $$\qq \qq\q \sum_{j=1}^k \frac{1}{p_j} 
 < \log\log \theta(p_k)
 +B_1 + p_k^{-0.5 +\varepsilon},
  \q \forall\ \ve > 0;
 \eqno(4.2)
 $$
 $$\qq \qq\q \sum_{j=1}^k \frac{1}{p_j} 
 > \log\log \theta(p_k)
 +B_1 - p_k^{-0.5 +\varepsilon},
  \q \forall\ \ve > 0;
 \eqno(4.3)
 $$
\vspace{-2ex}
 $$\qq\qq \q  \sum_{j=1}^k \frac{1}{p_j} 
 <\log\log \theta(p_k) + B_1 
 + \frac{A_0 }{\sqrt{p_k} \log p_k},\q A_0<\infty,
 \eqno(4.4)
 $$
 {\it where $B_1$ stands for the Meissel-Mertens constant:}
\vspace{-1ex}
 $$  \qq B_1:=\lim_{k\to\infty} \left(\sum_{j=1}^k \frac{1}{p_j} 
 - \log\log \theta(p_k)\right) 
=\gamma - \sum_p \left(\log \frac{p}{p-1} - \frac{1}{p}\ \right)
 $$
\vspace{-2ex}
 $$ 
=\gamma - \sum_p \sum_{k=2}^{\infty} \frac{1}{kp^k}
=0.261\ 497\dots
\eqno{(4.5)} 
 $$ 
 

For the proof one should use the Theorem and notice  that (cf. (1.2)) 
\vspace{-1ex}
 $$ \sum_{j=1}^k \left(\frac{1}{p_j}
  + \log\left(1-\frac{1}{p_j}\right)\right) =B_1 - \gamma  
  +O\left(\frac{1}{p_k}\right).
\eqno{(4.6)} 
 $$

3) Quite recently the author established (combining the method by Ingham [8, Sect. V. 10] with the properties of the so called locally G-maximal numbers, studied in [9, Sect 2]), that {\it if in (1.10)  $A_0<+\infty$,
 and thus {\bf RH} is true, then necessarily}
$$ 1.5-\ve\ <\ Q(x) \sqrt{x}\log x\ <\ 2.5+\ve, \ \forall\ x>X_{\ve}
\eqno{(4.6)} 
 $$ 
whence one may deduce  the relationship:

\vspace{-1ex}
$$ \max\{ \log G(N): {\bf gpf}(N)=p_k\}
= \gamma - \frac{a_k}{\sqrt{p_k}\log p_k}; 
$$
$$ 2\sqrt{2} - 2.5 -\ve < a_k < 2\sqrt{2} - 1.5 +\ve, \q \forall\ k>K_{\ve},
\eqno{(4.7)}  
$$ 
which quantitatively  refines the initial Ramanujan inequality (1.6).

\vspace{1ex}
  
These results will be presented in the next author's papers.

\vspace{1ex}

\vspace*{\fill}

\clearpage

\vspace{2ex}

\centerline{\bf LIST OF REFERENCES }

\vspace{3ex}

\noindent
[1] Mertens F. {\it $\ddot{U}ber$ einige asymptotische Gesetze der Zahlentheorie}.

 J. Reine Angew.  Math.,  {\bf 77}, 1874, pp. 289 -- 338.

\vspace{1ex}

\noindent
[2]\   Robin G. \ {\it Grandes valeurs de la fonction somme 
des diviseurs et hypoth\`ese }

{\it de  Riemann.} J. Math. Pures Appl. V. 63 (1984), pp. 187 -- 213.

\vspace{1ex}

\noindent
[3]\    Caveney G., Nicolas J.-L., Sondow J. \ {\it Robin's theorem,
 primes,  }

{\it and a new elementary reformulation of Riemann Hypothesis.}

INTEGERS 11 (2011),  \#A33, pp 1 -- 10.

\vspace{1ex}

\noindent
[4]\  Gronwall T.H. {\it Some asymptotic expressions 
in the theory of numbers}. 

Trans. Amer.  Math. Soc. V. 14 (1913), pp. 113 -- 122.

\vspace{1ex}
\noindent
[5]\   Ramanujan S. \ {\it Highly composite numbers, annotated 
and with a foreword} 

{\it by J.-L. Nicolas and G. Robin},
\ Ramanujan J. V. 1,  1997, pp. 119 -- 153.

\vspace{1ex}

\noindent
[6]\ Narkiewicz W. {\it The Development of Prime Number Theory}. 
 
{Springer-Verlag Berlin Heidelberg New York},  2000.

\vspace{1ex}

\noindent
[7]\   Nicolas J.-L. \  {\it Petites valeurs de la fonction d'Euler.} 

 Journal of Number Theory, vol. 17, no. 3, 1983, p. 375 -- 388.

\vspace{1ex}

\noindent
[8]\ Ingham A.\,E. {\it Distribution of Prime Numbers}. 
 
{\ Cambridge at the University Press},  1932.

\vspace{1ex}

\noindent
[9]\ Kalyabin G.\,A. {\it One-Step $G$-Unimprovable Numbers}. 
 
{\ https://arxiv.org/abs/1810.12585},  2018.

\vspace{1ex}

\vspace*{\fill}

\end{document}